\documentclass[a4paper]{amsart} 
\usepackage{amssymb,times, amscd,amsmath,amsthm, xypic}
\usepackage{geometry}
\title
{A remark on Yoneda's Lemma}

\author{Shoji Yokur$\mbox{a}^*$}
\thanks{
\\
(*) Partially supported by JSPS KAKENHI Grant Numbers 16H03936}

\date{}

\address{Department of Mathematics  and Computer Science, Graduate School of Science and Engineering, Kagoshima University, 1-21-35 Korimoto, Kagoshima, 890-0065, Japan
}

\email{yokura@sci.kagoshima-u.ac.jp}

\begin{document} 
\numberwithin{equation}{section}
\newtheorem{thm}[equation]{Theorem}
\newtheorem{pro}[equation]{Proposition}
\newtheorem{prob}[equation]{Problem}
\newtheorem{qu}[equation]{Question}
\newtheorem{cor}[equation]{Corollary}
\newtheorem{con}[equation]{Conjecture}
\newtheorem{ob}[equation]{Observation}
\newtheorem{lem}[equation]{Lemma}
\theoremstyle{definition}
\newtheorem{ex}[equation]{Example}
\newtheorem{defn}[equation]{Definition}
\newtheorem{rem}[equation]{Remark}
\renewcommand{\rmdefault}{ptm}
\def\alp{\alpha}
\def\be{\beta}
\def\jeden{1\hskip-3.5pt1}
\def\om{\omega}
\def\bigstar{\mathbf{\star}}
\def\ep{\epsilon}
\def\vep{\varepsilon}
\def\Om{\Omega}
\def\la{\lambda}
\def\La{\Lambda}
\def\si{\sigma}
\def\Si{\Sigma}
\def\Cal{\mathcal}
\def\ga{\gamma}
\def\Ga{\Gamma}
\def\de{\delta}
\def\De{\Delta}
\def\bF{\mathbb{F}}
\def\bH{\mathbb H}
\def\bPH{\mathbb {PH}}
\def \bB{\mathbb B}
\def \bA{\mathbb A}
\def \bOB{\mathbb {OB}}
\def \bM{\mathbb M}
\def \bOM{\mathbb {OM}}
\def \calB{\mathcal B}
\def \bK{\mathbb K}
\def \bG{\mathbf G}
\def \bL{\mathbf L}
\def\bN{\mathbb N}
\def\bR{\mathbb R}
\def\bP{\mathbb P}
\def\bZ{\mathbb Z}
\def\bC{\mathbb  C}
\def \bQ{\mathbb Q}
\def\op{\operatorname}


\maketitle

\begin{abstract} Yoneda'e Lemma is about the canonical isomorphism of all the natural transformations from a given representable covariant (contravariant, reps.) functor (from a locally small category to the category of sets) to  a covariant (contravariant, reps.) functor. In this note we point out that given any representable functor and any functor we have the canonical natural transformation from the given representable functor to the ``subset" functor of the given functor, ``collecting all the natural transformations".

\end{abstract}

\section {Yoneda's Lemma}

The well-known Yoneda's lemmas about representable functors are the following:
\begin{thm}Let $\mathcal C$ be a locally small category, i.e., $hom_{\mathcal C}(A,B)$ is a set, and let $\mathcal Set$ be the category of sets.

\begin{enumerate}
\item (the covariant case) Let $F_*:\mathcal C \to \mathcal Set$ be a coavariant functor. Let $h_A:= hom_{\mathcal C}(A, -)$ be a covariant hom-set functor $h_A: \mathcal C \to \mathcal Set$. Then the set of all the natural transformations from the hom-set covariant functor $h_A =hom_{\mathcal C}(A,-)$ to the covariant functor $F_*$ is isomorphic to the set $F_*(A)$:
$$\mathcal Natural (h_A, F_*) \cong F_*(A).$$

\item (the contravariant case) Let $F^*:\mathcal C \to \mathcal Set$ be a contravariant functor. Let $h^A:= hom_{\mathcal C}(-, A)$ be a contravariant hom-set functor $h^A: \mathcal C \to \mathcal Set$. Then the set of all the natural transformations from the hom-set contravariant functor $h^A =hom_{\mathcal C}(-, A)$ to the contravariant functor $F^*$ is isomorphic to the set $F^*(A)$:
$$\mathcal Natural (h^A, F^*) \cong F^*(A).$$
\end{enumerate}
\end{thm}

From now, for the sake of simplicity, we denote $h^A(X) =hom_{\mathcal C}(X, A)$ simply by $[X,A]$ and similarly $[A,X]$ for $h_A(X) = hom_{\mathcal C}(A,X)$.
 
The contravariant case of Yoneda's Lemma is proved by using the following commutative diagram: Let $\tau: [-, A] \to F_*(-)$ be a natural transformation:
\begin{equation}
\xymatrix
{
\op{id}_A\in [A, A] \ar[rr]^{\tau} \ar[d]_{f^*} && F^*(A) \ni \tau (\op{id}_A) \ar[d]^{f^*}\\
f \in [X, A] \ar[rr]_{\tau} && F^*(X)\ni \tau(f).
}
\end{equation}
Note that $f = f^*(\op{id}_A)= f \circ \op{id}_A$. Hence we have
\begin{align*}
\tau(f) & =  \tau (f^*(\op{id}_A)) \\
& = f^*(\tau (\op{id}_A)) \, \, \text {(by the naturality of $\tau$)}
\end{align*}
Thus the natural transformation $\tau: [-, A] \to F_*(-)$ is determined by the element $\tau(\op{id}_A) \in F^*(A)$. Conversely, given any element $\alp \in F^*(A)$ we can define the natural transformation $\tau_{\alp}:[-,A] \to F^*(-)$ by, for each object $X \in Obj(\mathcal C)$
$$\tau_{\alp}: [X, A] \to F^*(X) \quad \tau_{\alp}(f) = f^*\alp,$$
in which case $\tau_{\alp}:[A,A] \to F^*(A)$ satisfies that $\tau_{\alp}(\op{id}_A) = \op{id}_A^*(\alp)=\alp$.
The above isomorphism map is called the Yoneda map:
$$\mathcal Y: \mathcal Natural (h^A, F^*) \cong F^*(A) \quad \mathcal Y(\tau):=\tau(\op{id}_A),\text{or}$$
$$\mathcal Y: F^*(A) \cong  \mathcal Natural (h^A, F^*) \quad \mathcal Y(\alp):=\tau_{\alp}.$$

\section{a remark}
Let $G$ be a set and let $\mathcal Sub(G)$ be the set of all subsets of $G$.
Let $h^A(-)=[-,A]:C \to \mathcal Set$ and $F^*: \mathcal C \to \mathcal Set$ be as above. Then for each object $X\in Obj(\mathcal C)$ we have the following canonical map:
$$ \mathcal Im_{F^*}: h^A(X)=[X,A] \to \mathcal Sub (F^*(X))$$ defined by
$$\mathcal Im_{F^*}(f) := \op{Image} (f^*:F^*(A) \to F^*(X)) = f^*(F^*(A)) = \{f^*\alp \, | \, \alp \in F^*(A) \}.$$
The last two parts are written down for an emphasis. As observed in the above, $f^*\alp = \tau_{\alp}(f)$
which is the image of $f$ under the natural transformation $\tau_{\alp}$ corresponding to $\alp \in F^*(A)$.
In other word 

$\mathcal Im_{F^*}(f)$ \emph{is the set consisting of the images of $f$ by \underline {all} the natural transformations $\mathcal Natural (h^A, F^*)$}. 
For a morphism $g:X \to Y \in \mathcal C$, we have the following commutative diagram:

$$\xymatrix
{
h^A(Y)=[Y,A] \ar [rr]^{\mathcal Im_{F^*}} \ar [d]_{g^*} && \mathcal Sub (F^*(Y)) \ar [d]^{g^*} \\
h^A(X)=[X,A]  \ar[rr]_{\mathcal Im_{F^*}} && \mathcal Sub (F^*(X)) \\
}
$$
If we let $\mathcal Sub F^*: \mathcal C \to \mathcal Set$ be the ``subset" functor associated to the given functor $F^*:\mathcal C \to \mathcal Set$, defined by $\mathcal Sub F^*(X) := \mathcal Sub (F^*(X))$, we can consider $\mathcal Im_{F^*}(f)$ as \emph{a natural transformation}
$$\mathcal Im_{F^*}: h^A(-) \to \mathcal Sub F^* (-),$$
which \emph{sort of ``collects" all the natural transformation images}.

The upshot is the following. 
\begin{ob} Let $\mathcal C$ be a locally small category and $\mathcal Set$ be the category of sets. Let
$h^A(-)=[-, A]: \mathcal C \to \mathcal Set$ be a representable contravariant functor and $F^*:\mathcal C \to \mathcal Set$ be another contravariant functor. Then we have the following canonical natural transformation (`sort of ``collecting" or ``using" all the natural transformations from $h^A$ to $F^*$)
$$\mathcal Im_{F^*}:h^A(-) \to \mathcal Sub F^*(-),$$
where  for each object $X \in Obj(\mathcal C)$ we have $\mathcal Im_{F^*}(f) =f^*(F^*(A)) (\subset F^*(X) )$, which is the set consisting of the images of $f$ by \underline{all} the natural transformations from $h^A$ to $F^*$.
\end{ob}
The similar observation for the covariant case is made mutatis mutandis, so omitted.

\begin{rem} Depending on the situations, the target category of our contravariant functor can have more structures, e.g., groups, abelian groups, rings, commutative rings, etc.
\end{rem}

\section{Applications/Examples}
\begin{ex}[complex vector bundles and characteristic classes (e.g., see \cite{MS}, \cite{Ha}).] 
Let $\op{Vect}_n(X)$ be the set of isomorphism classes of complex vector bundles of rank $n$ and for $f:X \to Y$ the pullback map $f^*: \op{Vect}_n(Y) \to \op{Vect}_n(X)$ is define by $f^*[E]:=[f^*E]$. Thus the functor is a contravariant functor from the category of (paracompact) topological spaces to the category $\mathcal Set$ of sets.
Then we do know that
$$\op{Vect}_n(X) \cong [X, G_n(\mathbb C^{\infty})]$$
where $G_n(\mathbb C^{\infty})$ is the infinite Grassmann manifold of complex planes of dimension $n$, i.e., the classifying space of complex vector bundles of rank $n$.
This isomorphism is by the correspondence $[E] \longleftrightarrow [f_E]$, where $f_E:X \to G_n(\mathbb C^{\infty})$ is a classifying map of $E$, i.e., $E = f_E^* \gamma^n,$ 
where $\gamma^n$ is the universal complex vector bundle of rank $n$ over $G_n(\mathbb C^{\infty})$. 
Thus the functor $\op{Vect}_n(-)$ is a representable contravariant functor. Let $H^*(-; \mathbb Z)$ be the integral cohomology functor. Then by the observation in the previous section we have the following natural transformation:
$$ \op{Im}_{H^*}: \op{Vect}_n(-) \to \mathcal Sub H^*(-; \mathbb Z),$$
defined by for $\op{Vect}_n(X)$,
$ \op{Im}_{H^*}([E]):= \op{Im}\Bigl (f_E^*:H^*(G_n(\mathbb C^{\infty});\mathbb Z) \to H^*(X; \mathbb Z) \Bigr ).$
By the definition of characteristic classes, \emph{ $\op{Im}\Bigl (f_E^*:H^*(G_n(\mathbb C^{\infty});\mathbb Z) \to H^*(X; \mathbb Z) \Bigr )$ is nothing but the subring consisting of \underline{all} the characteristic classes of $E$}, which is $\mathbb Z [c_1(E), c_2(E), \cdots, c_n(E)]$. Let us denote this subring by $\op{Char}(E)$. By the definition for isomorphic two vector bundles $E$ and $E'$ we do have $Char(E)=Char(E')$.
One could define a very ``coarse" classification of vector bundles using $Char(E)$, i.e.,
$$E \cong_{coarse} E' \Longleftrightarrow Char(E) = Char(E').$$
For example, a line bundle $L$ and its inverse $L^{-1}$ satisfy that $L \cong_{coarse}L^{-1}$
because $c_1(L^{-1})=-c(L)$, which implies that $Char(L) = Char(L^{-1})$.
\end{ex}

\begin{rem} In the case of real vector bundles, the complex infinite Grassmann $G_n(\mathbb C^{\infty})$, the Chern class $c_i$ and the coefficient ring $\mathbb Z$ are respectively replaced by the real infinite Grassmann 
 $G_n(\mathbb R^{\infty})$, the Stiefell-Whitney class $w_i$ and the coefficient ring $\mathbb Z_2$.
\end{rem}

\begin{ex}[Ren\'e Thom's notion of \emph{dependence of cohomology classes}] Thom defined the following:
\begin{defn}[R. Thom] The cohomology class $\beta \in H^q(X; B)$ \emph{depends on} the cohomology class $\alpha \in H^p(X; A)$, where $A, B$ are coefficient rings, if, for all (perhaps infinite) polyhedra $Y$ and all maps $f: X \to Y$ such that $\alpha \in f^*(H^p(Y; A))$, we have $\beta \in f^*(H^q(Y;B))$.
\end{defn}

Fist we recall that the cohomology theory is a representable contravariant functor, indeed, representable by the Eilenberg-Maclane space, i.e., $H^j(X, R) \cong [X, K(R,j)]$ where $K(R,j)$ is the Eilenberg-Maclane space whose homotopy type is completely characterized by the homotopy groups $\pi_j(K(R,j)) = R$ and $\pi_i(K(R,i))= 0, i \not = j.$ Then by the Hurewicz Theorem we have $H_j(K(R,j); \mathbb Z)\cong \pi_j(K(R,j)) = R$ and $H_d(K(R,j)) = 0$ for $d<j$. Hence by the universal coefficinet theorem we have the isomorphism
$$\Phi:H^j(K(R,j);R) \cong Hom(H_j(K(R,j);\mathbb Z), R) \cong Hom (\pi_j(K(R,j)) , R) \cong Hom(R, R).$$
Let $u :=\Phi^{-1}(\op{id}_R)$ for the identity map $\op{id}_R: R \to R$. Then the isomorphism $\Theta: [X, K(R,j)] \cong H^j(X, R)$ is obtained by $\Theta ([f]):=f^*u$ where $f^*:H^j(K(R,j);R) \to H^j(X,R)$.
\begin{pro}[R. Thom] Let $\alpha \in H^p(X; A) \cong [X, K(A, p)]$ and let $f_{\alpha}:X \to K(A, p)$ be a map such that the homotopy class $[f_{\alpha}]$ corresponds to $\alpha$.
Then $\beta \in H^q(X,B)$ depends on $\alpha$ if and only if $\beta \in f^*_{\alpha}(H^q(K(A,p); B))$.
\end{pro}
In our set-up, we consider the representable contravariant functor $h^{K(A, p)}(-)$ and the contravariant cohomology functor $H^q(-; B)$. Then from the above section we have the following canonical natural transformation:
$$\mathcal Im_{H^q(-;B)}: [-, K(A, p)] \to \mathcal Sub H^q (-; B),$$
where for a topological space $X$ and for (the homotopy class) of $f_{\alp}:X \to K(A, p)$ corresponding to the cohomology class $\alp \in H^p(X; A) (\cong [X, K(A,p)])$ we have
$$\mathcal Im_{H^q(-;B)}(f_{\alp}) = f_{\alp}^*(H^q(K(A,p); B)),$$
which is, due to the above proposition of Thom, \emph{nothing but the subgroup of all the cohomology classes $\beta \in H^q(X;B)$ depending on the cohomology class $\alpha$}, and also by our observation above it is \emph{the subgroup consisting of the image of $f_{\alp}$ by all the natural transformations from the representable functor $[-, K(A, p)]$, in other words the cohomology functor $H^p(-;A)$ to the cohomology functor $H^q(-;B)$}.
\end{ex}

One can consider some other reasonable or interesting pairs $(h^A(-), F^*(-))$ of representable contravariant functors $h^A(-)$ and contravariant functors $F^*(-)$. In a different paper we want to study such things, e.g. $K$-theory and the cohomology theory.

\section{one more remark: poset-stratified space structures of representable functors}

The subset functor $\mathcal Sub F^*(-)$ has in fact another simple structure of partial order set due to the set-theoretic inclusion. If we consider the above natural transformation
$$ \mathcal Im_{F^*}: h^A(X)=[-,A] \to \mathcal Sub F^*(-)$$
and for an object $X \in Obj(\mathcal C)$ we have
$$ \mathcal Im_{F^*}: h^A(X)=[X,A] \to \mathcal Sub (F^*(X))$$
defined by $\mathcal Im_{F^*}(f) = f^*(F^*(A))$. For the following commutative diagram of morphisms in $\mathcal C$
$$
\xymatrix
{ X \ar[rr]^f \ar[dr]_g && A \\
& A \ar[ur]_t
}
$$
under the contravaiant functor $F^*$, we have that $f^*= g^* \circ t^*$:
$$
\xymatrix
{ F^*(X) && F^*(A) \ar[ll]_{f^*} \ar[dl]^{t^*}\\
& F^*(A) \ar[ul]^{g^*}
}
$$
Hence we have that 
$$\mathcal Im_{F^*}(f) = f^*(F^*(A)) \subseteq \mathcal Im_{F^*}(g) = g^*(F^*(A)).$$
If we define the order $f \leqq_L g$ for $f, g \in [X, A]$ by the above commutative diagram $
\xymatrix
{ X \ar[rr]^f \ar[dr]_g && A \\
& A \ar[ur]_t
}
$, i.e., by the condition that $\exists t \in [A,A]$ such that $f = t \circ g$, then this order $\leqq_L$  is a preorder, i.e., it is reflexive and transitive, but not necessarily anti-symmetric. With this order we get a preordered set $([X, A], \leqq_L)$. A preordered set is called \emph{a proset} and the category of presets and monotone (order-preserving) maps is denoted by $\mathcal Proset$. Namely, by this preorder, the representable contravariant functor $h^A: \mathcal C \to \mathcal Set$ becomes a representable contravariant functor to the category of prosets: $h^A: \mathcal C \to \mathcal Proset$. The contravariant functor $F^*(-)$ also gives rise to the associated contravariant functor $\mathcal Sub F^*: \mathcal C \to \mathcal Proset$, to be more precise, the order of inclusion is a partial order, thus in the case of $\mathcal Sub F^*$ the target category id the category $\mathcal Poset$ of posets (partially ordered sets) and monotone (order-preserving) maps. With these preorders the above natural transformation $ \mathcal Im_{F^*}: h^A(X)=[-,A] \to \mathcal Sub F^*(-)$ is a natural transformation $ \mathcal Im_{F^*}: (h^A(X)=[-,A], \leqq_L) \to (\mathcal Sub (F^*(-)), \leqq))$.
Then for each object $X \in Obj(\mathcal C)$ we have the map
$$ \mathcal Im_{F^*}: (h^A(X)=[X,A], \leqq_L) \to (\mathcal Sub (F^*(X)), \leqq))$$
which is a monotone map from a proset to a poset. If we consider the Alexandroff topologies \cite{A}(also see \cite{Ar, AFT, Sp}) for a proset, this monotone map becomes a continuous map from the proset considered as a topological space with the Alexandroff topology to a poset considered as a topological space with the Alexandroff topology. This is nothing but the so-called \emph{a poset-stratified space} \cite{Lurie}(e.g., also see \cite{AFT}). From this viewpoint in \cite{YY2} (cf. \cite{YY}) we consider poset-stratified space structures of the homotopy set of continuous maps of topological spaces and in \cite{Yo} for a general locally small category.

\begin{rem} To get a poset-stratified space structure of $ \mathcal Im_{F^*}: (h^A(X), \leqq_L) \to (\mathcal Sub (F^*(X)), \leqq))$ we appeal to another contravariant functor $F^*(-)$. But, in order to get such a poset-stratified space structure we do not need such a functor. In general, if we have a proset $(P, \leqq)$, then we consider the equivalence relation $a \sim b \Leftrightarrow a \leqq b, b \leqq a$ on $P$ and consider the set $P_{\sim}$ of the equivalence relations and we define the order $[a] \leqq' [b] \leftrightarrow a \leqq b$, which is a partial order and the canonical map $\pi: (P, \leqq) \to (P_{\sim}, \leqq')$ is a monotone map, thus we get a poset-stratified space considering the associated Alexandroff topologies (e.g., \cite{Yo}). This construction give a kind of universal poset-stratified space structure $\pi: ([X,A], \leqq_L) \to ([X,A]_{\sim}, \leqq_L')$ and the above $ \mathcal Im_{F^*}: (h^A(X)=[X,A], \leqq_L) \to (\mathcal Sub (F^*(X)), \leqq))$ involving another contravariant functor $F^*$ give a more geometric one, so to speak.
\end{rem}
\begin{rem} For the topological homotopy category $h\mathcal Top$ the above relation $f = t \circ g$ becomes
$f \sim t \circ g$. This relation was considered in a different context by K. Borsuk \cite{Bor1, Bor2} and generalized by P. Hilton \cite{Hil1} (cf. \cite{Hil2, Hil3}).

\end{rem}

\end{document}